\documentclass[12pt,psamsfonts]{amsart}

\newtheorem{anyprop}{Anyprop}[section]

\newtheorem{lemma}[anyprop]{Lemma}
\newtheorem{proposition}[anyprop]{Proposition}
\newtheorem{corollary}[anyprop]{Corollary}

\theoremstyle{definition}

\newtheorem{example}[anyprop]{Example}

\newtheorem{remark}[anyprop]{Remark}

\newcommand{\NN}{\mathbb{N}}
\newcommand{\ZZ}{\mathbb{Z}}

\newcommand  {\shE}     {\mathcal{E}}

\newcommand  {\foA}     {\mathfrak{A}}

\newcommand  {\Char}    {\operatorname{char}}

\newcommand  {\lra}     {\longrightarrow}

\newcommand  {\modu}     {\operatorname{mod}}

\renewcommand{\O}       {\mathcal{O}}

\newcommand  {\ord}     {\operatorname{ord}}

\newcommand  {\Proj}    {\operatorname{Proj}}

\newcommand  {\ra}      {\rightarrow}

\newcommand  {\rk}    {\operatorname{rk}}

\newcommand  {\Spec}    {\operatorname{Spec}}

\newcommand  {\Syz}     {\operatorname{Syz}}

\newcommand{\comdots}{ , \ldots , }
\newcommand{\komdots}{ , \ldots , }
\newcommand{\plusdots}{ + \ldots + }

\theoremstyle{remark}

\numberwithin{equation}{section}

\usepackage{amscd}
\usepackage{amssymb}

\def\mydate{\number\day\space\ifcase\month \or January\or February\or March\or April\or May\or
June\or July\or August\or September\or October\or November\or
December\fi \space\number\year}

\newcommand{\ellkk}{k}

\newcommand{\elll}{\ell}

\newcommand{\soger}{h}

\newcommand{\rangem}{M}

\newcommand{\numb}{\sharp}

\begin{document}

\title[Miyaoka problem]
{On a problem of Miyaoka}

\author[Holger Brenner]{Holger Brenner}
\address{Department of Pure Mathematics, University of Sheffield,
  Hicks Building, Hounsfield Road, Sheffield S3 7RH, United Kingdom}
\email{H.Brenner@sheffield.ac.uk}





\begin{abstract}
We give an example of a vector bundle $\shE$ on a relative curve $C
\ra \Spec \ZZ$ such that the restriction to the generic fiber in
characteristic zero is semistable but such that the restriction to
positive characteristic $p$ is not strongly semistable for
infinitely many prime numbers $p$. Moreover, under the hypothesis
that there exist infinitely many Sophie Germain primes, there are
also examples such that the density of primes with non strongly
semistable reduction is arbitrarily high.
\end{abstract}

\maketitle


\section*{Introduction}

In this paper we deal with the following problem of Miyaoka \cite[Problem 5.4]{miyaokachern}:

``Let $C$ be an irreducible smooth projective curve over a noetherian integral domain $R$ of
characteristic $0$.
Assume that a locally free sheaf $\shE$ on $C$ is $\foA$-semistable on the generic fiber $C_*$.
Let $S$ be the set of primes of char. $>0$ on $\Spec R$
such that $\shE$ is strongly semistable.
Is $S$ a dense subset of $\Spec R$?''

Here $\shE$ is called strongly semistable if every Frobenius pull-back
of $\shE$ stays semistable.
Since semistability is an open property and since
semistable bundles on the projective line and on an elliptic curve
are strongly semistable, this problem has a positive answer for genus $g=0,1$.

Shepherd-Barron `rephrases' the question asking
whether ``is it true
that the set ... of primes $p$ modulo which $\shE$ is not strongly semistable (...)
is finite, or at least of density $0$?'' \cite{shepherdbarronsemistability}.
He considers also higher dimensional varieties $V$,
and one of his main results is that for $\dim V \geq 2$ and $\rk (\shE)=2$
the set $\Sigma$ of prime numbers with non strongly semistable reduction
is finite under the condition that either
the Picard number of $V$ is one, or that the canonical bundle $K_V$ is numerically trivial,
or that the variety $V$ is algebraically simply connected \cite[Corollary 6]{shepherdbarronsemistability}.

Coming back to curves of genus $\geq 2$, say over $\Spec \ZZ$, nothing seems to be known about the following questions on the set $S$ of prime numbers with strongly semistable
reduction: Does $S$ contain almost all prime numbers (as in the results of Shepherd-Barron)?
Is $S$ always an infinite set? Is it even possible that $S$ is empty?
Can we say anything about the density of $S$ in the sense of analytic number theory?

In this paper we give examples of vector bundles of rank two, which
are semistable in characteristic zero, but not strongly semistable
for infinitely many prime numbers. We also provide examples where
the density of primes with non strongly semistable reduction is very
high, and in fact arbitrarily high under the hypothesis that there
exists infinitely many Sophie Germain prime numbers.

The example is just the syzygy bundle $\Syz(X^2,Y^2,Z^2)$
on the plane projective curve given by
$Z^d=X^d+Y^d$ for $d \geq 5$.
These bundles are semistable in characteristic zero.
The point is that in positive characteristic $p$ fulfilling certain congruence conditions
modulo $d$,
some Frobenius pull-backs of these bundles have global sections which contradict the strong semistability.
It is also possible that in these examples the reduction is not strongly semistable
for all prime numbers, but this we do not know.

This type of examples is motivated by the theory of tight closure.
It was already used in \cite{brennerstronglysemistable} to show that
there is no Bogomolov type restriction theorem for strong
semistability in positive characteristic. We will come back to the
impact of these examples on tight closure and on Hilbert-Kunz theory
somewhere else. I thank Neil Dummigan (University of Sheffield) for
a useful remark concerning Sophie Germain primes.

\section{Main results}

In the following we will consider syzygy bundles of rank two on a smooth
projective curve $C= \Proj A$ over a field $K$, where $A$ is a two-dimensional normal standard-graded $K$-domain.
Such a bundle is given by three homogeneous, $A_+$-primary elements $f_1,f_2,f_3 \in A$
of degree $d_1,d_2,d_3$ by the short exact sequence
$$ 0 \! \ra \! \Syz(f_1,f_2,f_3)(m) \ra \O(m-d_1) \oplus \O(m-d_2) \oplus \O(m-d_3)
\stackrel{f_1,f_2,f_3}{\ra} \O(m) \! \ra \! 0 .$$
A global section of
$\Syz(f_1,f_2,f_3)(m)$ is a triple $(s_1,s_2,s_3)$ of homogeneous elements
such that $\deg(s_i) + d_i =m$, $i=1,2,3$, and $s_1f_1+s_2f_2+s_3f_3=0$.
We call $m$ the total degree of the syzygy $(s_1,s_2,s_3)$.
The degree of such a syzygy bundle is by the short exact sequence
$ \deg (\Syz(f_1,f_2,f_3)(m)) =(2m -d_1-d_2-d_3) \deg (C)$.
If $m$ is such that this degree is negative and such that there exists global non-trivial sections,
then this bundle is not semistable.

It is in general not easy to control the global syzygies; in the following lemma however we take advantage of the existence of a noetherian normalization of a very special type.
Let $\delta (f_1 \komdots f_n)$ denote the minimal total degree of a non-trivial syzygy for
$f_1 \komdots f_n$.

\begin{lemma}
\label{syzygylemma}
Let $P(X,Y)$ denote a homogeneous polynomial in $K[X,Y]$ of degree $d$ and
consider the projective curve $C$ given by $Z^d-P(X,Y)$,
so that $C=\Proj K[X,Y,Z]/(Z^d-P(X,Y))$. Suppose that $C$ is smooth.
Let $a_1,a_2,a_3 \in \NN$
and consider the syzygy bundle
$\Syz(X^{a_1},Y^{a_2}, Z^{a_3})(m)$ on $C$.
Write $a_3=dk+t$ where $0 \leq t < d$.
Then
$$ \delta( X^{a_1},Y^{a_2}, Z^{a_3}) = \min \{   \delta( X^{a_1},Y^{a_2}, P(X,Y)^k ) + t
, \, \delta( X^{a_1},Y^{a_2}, P(X,Y)^{k+1} ) \} \, .$$
\end{lemma}
\proof A global section of $\Syz(X^{a_1},Y^{a_2}, Z^{a_3})(m)$ is
the same as homogeneous polynomials $F,G,H \in
K[X,Y,Z]/(Z^d-P(X,Y))$ such that $FX^{a_1}+ GY^{a_2}+HZ^{a_3}=0$ and
$\deg (F)+a_1=\deg (G)+a_2 = \deg (H)+a_3=m$. We may write $F=F_0+
F_1Z+F_2Z^2 \plusdots F_{d-1}Z^{d-1}$ with $F_i \in K[X,Y]$ and
similarly for $G$ and $H$. We have $Z^{a_3}= Z^{dk+t} = (Z^d)^kZ^t
=P(X,Y)^k Z^t$. The polynomials $(F,G,H)$ (fulfilling the degree
condition) are a syzygy if and only if for $i=0 \komdots d-1$ we
have
$$F_iZ^{i}X^{a_1} + G_iZ^{i} Y^{a_2}  +H_j Z^{j} Z^{a_3} =0,\,
\mbox{ where } j =i-t \modu d \, .$$ Let now $(F,G,H)$ denote a
non-trivial syzygy of minimal degree for $ X^{a_1}$, $Y^{a_2}$,
$Z^{a_3}$ and let $i$ denote the minimal number such that $F_i \neq
0$ or $G_i \neq 0$. Since the degree is minimal we may assume by
dividing through $Z^{\min (i,j)}$ that either $i=0$ or $j=0$, which
means $i=t$.

In the first case we can read the zero-component of the syzygy directly as a non-trivial
syzygy for  $X^{a_1},Y^{a_2},P(X,Y)^{k+1}$.
In the second case the $i=t$-th component of the syzygy is
$$F_tZ^{t}X^{a_1} + G_tZ^{t} Y^{a_2} +  H_{0} Z^{a_3} =0 \, .$$
Replacing  $Z^{a_3}$ through $P(X,Y)^k Z^t$  and dividing through $Z^t$ we get a non-trivial
syzygy for $X^{a_1},Y^{a_2},P(X,Y)^{k}$ of the same degree $-t$.

Suppose that we have a non-trivial syzygy
$ \tilde{F} X^{a_1}+ \tilde{G}  Y^{a_2}+ \tilde{H} P(X,Y)^{k+1}=0 $ in $K[X,Y]$.
Then $F=F_0 = \tilde{F}$, $G=G_0= \tilde{G}$ and $H=H_{d-t}Z^{d-t} = \tilde{H} Z^{d-t}$
gives a syzygy for $X^{a_1}$, $Y^{a_2}$, $Z^{a_3}$
of the same degree.

Suppose that we have a non-trivial syzygy
$\tilde{F} X^{a_1}+ \tilde{G}  Y^{a_2}+ \tilde{H} P(X,Y)^{k}
=0 $ in $K[X,Y]$.
Multiplying with $Z^t$ we see that
$F=F_{t}Z^t= \tilde{F} Z^{t}$, $ G=G_{t}Z^t= \tilde{G} Z^{t}$
and $H= \tilde{H}$
gives a syzygy for $X^{a_1}$, $Y^{a_2}$, $Z^{a_3}$
of the same degree $+ t$.
\qed

\begin{proposition}
\label{numcrit}
Let $d$ and $b$ denote natural numbers,
write $b =d \ellkk + t$ with $0 \leq  t < d$.
Suppose that $\ellkk $ is even and that $t  > 2d/3 $.
Then the syzygy bundle
$\Syz(X^{b}, Y^{b},Z^{b})$ is not semistable on the Fermat curve given by $Z^d =X^d +Y^d$.
\end{proposition}
\proof We will look for syzygies for $X^b$, $Y^b$ and
$(X^d+Y^d)^{\ellkk +1}$ of total degree $d(\ellkk+1) +d \lfloor
\ellkk /2 \rfloor =d( \ellkk +1 + \lfloor \ellkk /2 \rfloor)$, which
yields syzygies for $X^b,Y^b,Z^b$ of the same degree by Lemma
\ref{syzygylemma}. To find such syzygies we have to look for
multiples $H (X^d+Y^d)^{\ellkk +1} \in (X^{b},Y^{b} )$, where $\deg
(H)= d\lfloor \ellkk/2 \rfloor $. We consider for $H$ only monomials
in $X^{d }$ and $Y^{d}$, so these are the $\lfloor \ellkk/2 \rfloor
+1$ monomials
$$X^{d  \lfloor \ellkk /2 \rfloor}, X^{d(\lfloor \ellkk/2 \rfloor-1)}Y^d, X^{d(\lfloor \ellkk/2 \rfloor-2)}Y^{d2}
\komdots  Y^{d\lfloor \ellkk/2 \rfloor} \, .$$
The resulting monomials in the products, which do not belong to the ideal $(X^{b},Y^{b})$, have the form
$X^{d u}Y^{dv}$ with $du+dv=d( \ellkk + 1 +  \lfloor \ellkk/2 \rfloor   )$
and $u, v < d(\ellkk +1)$.
So these are the monomials
$$  X^{d( \lfloor \ellkk/2 \rfloor +1 )} Y^{d \ellkk },
X^{d( \lfloor \ellkk/2 \rfloor+2)}Y^{d(\ellkk -1)} \komdots  X^{d \ellkk }Y^{d(\lfloor \ellkk/2 \rfloor  +1)}\, . $$
These are $ \ellkk -(\lfloor \ellkk/2 \rfloor +1)   +1 = \ellkk - \lfloor \ellkk/2 \rfloor $ monomials.
Since $\ellkk$ is supposed to be even, we have
$ \lfloor \ellkk/2 \rfloor +1  >   \ellkk - \lfloor \ellkk/2 \rfloor$
and therefore we must have a non-trivial linear relation
$$\sum_{i+j= \lfloor \ellkk /2 \rfloor } \lambda_{ij} X^{di} Y^{dj} (X^d+Y^d)^{\ellkk +1}  = 0 , \,$$
modulo $(X^{b},Y^{b})$. The total degree of this non-trivial global syzygy is
$d( \ellkk +1 + \lfloor \ellkk /2 \rfloor)$
and the degree of the bundle is
$$
\deg \big( \Syz(X^{b}, Y^{b},Z^{b}) ( d( \ellkk +1 + \lfloor \ellkk
/2 \rfloor) ) \big) =
 \big(2 d( \ellkk +1 + \lfloor \ellkk /2\rfloor) -3b  \big)\deg (C)
$$
Due to our assumptions we have
$$2 d ( \ellkk +1 + \lfloor \ellkk/2 \rfloor) = 3d \ellkk +2d < 3d \ellkk +3t =3b \, ,$$
hence the degree of the bundle is negative, but it has a non-trivial section,
so it is not semistable.
\qed

\begin{corollary}
\label{numcrittwo}
Let $d$ and $q$ denote natural numbers,
write $q =d \elll +s$ with $0 \leq s < d$.
Suppose that $2s < d < 3s $.
Then the syzygy bundle
$\Syz(X^{2q}, Y^{2q},Z^{2q})$ is not semistable on the Fermat curve given by $Z^d =X^d +Y^d$.
\end{corollary}
\proof We apply Proposition \ref{numcrit} to $b=2q = d(2 \elll )
+2s= d \ellkk + t$. Note that $t < d$ and $2d < 6s=3t$. \qed

\begin{corollary}
\label{corollarynotstrong}
Consider the syzygy bundle $\shE=\Syz(X^2,Y^2,Z^2)$ on the Fermat curve
$C_K=\Proj K[X,Y,Z]/(X^d+Y^d-Z^d)$, $K$ a field.
Then $\shE_K$ is semistable in characteristic zero for $d \geq 5$,
but $\shE_K$ is not strongly semistable in positive characteristic $p =r \modu d $
such that some power $s =r^{e}$ fulfills $2s < d < 3s$.
In particular, for prime numbers $d \geq 5 $,
$\shE_K$ is not strongly semistable for infinitely many prime numbers $p$.
\end{corollary}
\proof Suppose first that $K$ has characteristic zero. Then $\shE_K$
is semistable due to \cite[Proposition
6.2]{brennercomputationtight}; this follows for $d \geq 7$ also from
the restriction theorem of Bogomolov, see \cite[Theorem
7.3.5]{huybrechtslehn}, since the bundle is clearly stable on the
projective plane.

Suppose now that $K$ has positive characteristic $p$ fulfilling the
assumption. Then we look at $q=p^{e} $ so that $q = d \elll  +s $
with $2s < d < 3s$, and Corollary \ref{numcrittwo} yields that
$\Syz(X^{2q},Y^{2q},Z^{2q})$ is not semistable. Since this bundle is
the pull-back under the $e$-th Frobenius of
$\Syz(X^{2},Y^{2},Z^{2})$, as follows from the short exact sequence
mentioned at the beginning of this section, we infer that
$\Syz(X^{2},Y^{2},Z^{2})$ is not strongly semistable. For prime
numbers $d \geq 5$ there are natural numbers $s$ such that $2s < d <
3s$ and such that $s$ is coprime to $d$. Due to the theorem of
Dirichlet about primes in an arithmetic progression \cite[Chapitre
VI, \S4, Th\'{e}or\`{e}me and Corollaire]{serrearithmetic}, there
exists infinitely many prime numbers $p$ with remainder $=s \modu
d$. \qed

\begin{remark}
The condition in Corollary \ref{corollarynotstrong} that $d$ is a
prime number is necessary, since for $d=6$ and $d=10$ there does not
exist such a coprime reminder $s$ in the range $d/3 < s < d/2$. Are
these the only exceptions?
\end{remark}

\begin{remark}
\label{remainderone} For $p=1 \modu d$ we have $q= 1 \modu d$ for
all $q=p^{e}$ and so Corollary \ref{corollarynotstrong} does not
apply. It is open whether for these prime numbers the bundle is
strongly semistable or not.

There is however some reason to believe that also in this case the
bundle is not strongly semistable. Suppose that we have $q= d \elll
+1$, hence $2q= d(2 \elll) + 2$. We look for syzygies for
$X^{2q},Y^{2q}, (X^d+Y^d)^{2 \elll}$ of total degree $ d(3 \elll)$.
This yields by Lemma \ref{syzygylemma} global syzygies for
$X^{2q},Y^{2q},Z^{2q}$ of degree $ d (3 \elll) +2$, which
contradicts the semistability, since $ 2 (d (3 \elll) +2) - 3( d(2
\elll)+2 ) =-2 $.

To find such syzygies we have to multiply $(X^d+Y^d)^{2 \elll}$ by
the $ \elll+1$ monomials $X^{d \elll}, X^{d( \elll-1)}Y^d \komdots
Y^{d \elll}$. The products are polynomials in the $ \elll +1$
monomials (modulo $(X^{2q},Y^{2q})$) $X^{d \elll}Y^{d2 \elll}, X^{d(
\elll +1)}Y^{d(2 \elll-1)} \komdots X^{d2\elll}Y^{d \elll} $. The
existence of such a non-trivial syzygy is equivalent to the property
that the determinant of the corresponding $(\elll+1) \times
(\elll+1)$-matrix is $0 \modu p$. Since all the prime powers
$q=p^{e}$ yield infinitely many such situations, it seems likely
that for some $e$ the determinant is zero.
\end{remark}

\section{The example on the Fermat quintic}

In this section we take a closer look at the example $\shE=
\Syz(X^2,Y^2,Z^2)$ on the Fermat quintic $Z^5=X^5+Y^5$ for various
characteristics. Then $\shE_K$ is semistable in characteristic zero,
but $\shE_K$ is not strongly semistable in characteristic $p=2$ or
$=3$ $\modu 5$. In characteristic $p=1$ or $p=4$ $\modu 5$ this is
not known.

\begin{corollary}
\label{corollarynotstrongfive}
Consider the syzygy bundle $\shE=\Syz(X^2,Y^2,Z^2)$ on the Fermat quintic
$C_K=\Proj K[X,Y,Z]/(X^5+Y^5-Z^5)$, $K$ a field.
Then $\shE_K$ is semistable in characteristic zero,
but $\shE_K$ is not strongly semistable in characteristic $p=2$ or $=3$ $\modu 5$.
For $p=2 \modu 5$ the first Frobenius pull-back of $\shE_K$ is not semistable,
and for $p=3 \modu 5$ the third Frobenius pull-back of $\shE_K$ is not semistable.
\end{corollary}
\proof For $p=3 \modu 5$ we have $q=p^3=2 \modu 5$, and for $p=2
\modu 5$ we take $q=p$. So in both cases we get a situation treated
in Corollary \ref{corollarynotstrong}, hence $\Syz(X^2,Y^2,Z^2)$ is
not strongly semistable.
\qed

\begin{remark}
In the case $p= 3 \modu 5$, Corollary \ref{corollarynotstrongfive}
shows that the third Frobenius pull-back of the syzygy bundle is not
semistable anymore. We show now that already the first Frobenius
pull-back is not semistable. Write $p=5u+3$ so that $2p= 5(2u+1) +1$
(and $k=2u+1$, $t=1$ in the notation of Lemma \ref{syzygylemma}). We
consider the syzygies for
$$  X^{5(2u+1) +1},\, Y^{5(2u+1) +1},\,  (X^5+Y^5)^{2u+1} \, .$$
We multiply the last term by the $u+1$ monomials
$$XY(X^{5u}Y^{0}), \,XY(X^{5(u-1)}Y^{5}) \comdots XY(X^0Y^{5u}) \, .$$
The resulting polynomials are modulo the first two terms expressible
in the monomials $X^{5i+1}Y^{5(3u+1-i)+1}$, where $i \leq 2u$ and
$u+1 \leq i$, so these are only $u$ many. Hence there exists a
syzygy of these polynomials of degree $5(3u+1) +2$ and therefore
there exists a global non-trivial syzygy for $X^{2p},Y^{2p},Z^{2p}$
of degree $5(3u+1)+3$ by Lemma \ref{syzygylemma}. This contradicts
semistability, since $2(5(3u+1)+3)-3(5(2u+1)+1)=-2$.

For example, for $p=3$, we find for $X^6,Y^6,(X^5+Y^5)$ the syzygy
$(-Y,-X,XY)$ of total degree $7$, which yields the syzygy
$(-YZ,-XZ,XY)$ for $X^6,Y^6,Z^6$ of total degree $8$ on the Fermat
quintic.
\end{remark}

\begin{example}
We consider the bundle $\Syz(X^2,Y^2,Z^2)$ on the Fermat quintic for $q=p=7$.
Then $Z^{14} =Z^{10}Z^4=(X^5+Y^5)^2 Z^4$
and we look for syzygies in $K[X,Y]$ for
$$X^{14},Y^{14}, (X^5+Y^5)^3 \, .$$
We multiply $(X^5+Y^5)^3$ by the monomials
$X^5$ and $Y^5$. The only monomial in the products which remains modulo $(X^{14},Y^{14})$
is $X^{10}Y^{10}$.
Therefore we must have a non-trivial syzygy of total degree $20$,
and indeed we have
$$ -( X^6+2XY^5)X^{14} + ( 2X^5Y+ Y^{6}) Y^{14}+ (X^5-Y^5)(X^5+Y^5)^3 =0  \, .$$
Going back to our original setting on the Fermat curve we get the syzygy
$$ -( X^6+2XY^5)X^{14} + ( 2X^5Y+ Y^{6}) Y^{14}+ (X^5-Y^5)Z  Z^{14} =0 \, .$$
This shows that $\Syz(X^{14},Y^{14},Z^{14})(20)$ has a non-trivial global section,
but its degree is $(2 \cdot 20 -3 \cdot 14) \deg(C) = -2 \deg (C)$ negative.
So $\Syz(X^2,Y^2,Z^2)$ is not strongly semistable for $p=7$.
It is easy to see that the syzygy $(- X^6-2XY^5, + 2X^5Y + Y^{6}, (X^5-Y^5)Z)$
does not have a common zero on the curve,
hence we get the short exact sequence
$$ 0 \lra \O  \lra \Syz(X^{14},Y^{14},Z^{14})(20)\lra \O(-2) \lra 0 \, ,$$
which is the Harder-Narasimhan filtration.
\end{example}

\begin{example}
We consider the example for $p=11$, so the remainder $\modu 5$ is
$1$ and we cannot expect a syzygy for $X^{22},Y^{22},Z^{22}$
contradicting the semistability. We have $Z^{22}=(X^5+Y^5)^4Z^2$ and
we look first for syzygies for $X^{22},Y^{22}, (X^5+Y^5)^4 $. We
have $(X^5+Y^5)^4 =X^{20} + 4X^{15}Y^{5} +6 X^{10}Y^{10} + 4
X^{5}Y^{15} + Y^{20}$ and multiplication by $X^{10}$, $X^5Y^5$,
$Y^{10}$ yields modulo $(X^{22},Y^{22})$ the three polynomials
$$ 6X^{20}Y^{10} +4 X^{15}Y^{15} +X^{10} Y^{20} \, ,$$
$$ 4X^{20}Y^{10} + 6 X^{15}Y^{15} + 4X^{10} Y^{20} \, ,$$
$$ X^{20}Y^{10} +4 X^{15}Y^{15} + 6X^{10} Y^{20} \, .$$
The determinant of the corresponding matrix is $50=6 \modu 11$
and so these polynomials are linearly independent.

We look now at syzygies for $X^{22},Y^{22}, (X^5+Y^5)^5 $. If we
consider only powers of $5$, we multiply only by $X^{5}$ and $Y^5$,
which yields modulo $(X^{22}, Y^{22})$ the monomials $ 10
X^{20}Y^{10} +10 X^{15}Y^{15} +5 X^{10}Y^{20}$ and $ 5 X^{20}Y^{10}
+10 X^{15}Y^{15} +10 X^{10}Y^{20}$, which are again linearly
independent.
\end{example}

\begin{remark}
For Fermat curves of degree $d<5$, the situation is as follows:
for $d=1$ the restriction of $\Syz(X^2,Y^2,Z^2)(3) = \Syz(X^2,Y^2,X^2+2XY+Y^2)(3)$
is $ \O \oplus \O$,
hence (strongly) semistable (characteristic $\neq 2$).
For $d=2$ the restriction of $\Syz(X^2,Y^2,Z^2)(2 )= \Syz(X^2,Y^2,X^2+Y^2)(2)$
has a non-trivial section, hence $\Syz(X^2,Y^2,Z^2)(3) \cong \O(-1) \oplus \O(1)$,
so this is not semistable in any characteristic.
For $d=3$ the Fermat equation $Z^3=X^3+Y^3$ yields at once
a global section $\O \ra \Syz(X^2,Y^2,Z^2)(3)$ without a zero.
This shows that the bundle is strongly semistable, but not stable,
independent of the characteristic.
For $d=4$ we have shown in \cite[Example 7.4]{brennercomputationtight}
that for $\Char (K) \neq 2$ the restriction is strongly semistable.
\end{remark}

\section{Using Sophie Germain primes}

Do there exist examples of vector bundles which are semistable in characteristic zero
and where the density of prime numbers for which the bundle
is not strongly semistable is arbitrarily high?
The density of prime numbers might be the analytic (or Dirichlet) density or the natural density
(see \cite[Chapitre 5, \S 4.1 and \S 4.5]{serrearithmetic}).
Since we will only use the fact
that the set of prime numbers $p$
such that $p=r \modu d$, ($r, d$ coprime),
has the density $1/ \varphi(d)$,
we will not say much about this point.

If we want to attack this question with the method of the first
section, we need to know for which and for how many remainders $r
\in \ZZ_d^\times$ there exists a power
$$r^{e} \in \rangem = \{s: d/3 < s < d/2\} \, .$$
For $r=1$ or $=-1 \modu d$ this is not possible; on the other hand,
it is always true for primitive elements if $\rangem $ is not empty.
For a remainder $r$ there exists some power inside $\rangem$ if and
only if the (multiplicative) group generated by $r$ intersects
$\rangem$. Therefore we only have to count the number of generators
of all the subgroups of $\ZZ_d^\times$ which contain an element of
$\rangem$, hence
$$ \numb \{r \in \ZZ_d^\times :\, \exists e \, \, r^{e} \in \rangem \}
= \sum_{H \subseteq \ZZ_d^\times, \, H \cap \rangem \neq \emptyset} \varphi(\ord (H)) \, ,
$$
where $\varphi$ denotes the Euler $\varphi$-function.
Good candidates to obtain here a big cardinality are degrees $d$ of type $d=2 \soger+1$, where both $d$ and $\soger$ are prime.
The numbers $\soger$ with this property are called Sophie Germain primes.
It is still not known whether there exist infinitely many such numbers.

\begin{proposition}
\label{germainfact} Suppose that $\soger > 5$ is a Sophie Germain
prime, set $d=2\soger +1$. Then the primes for which
$\Syz(X^2,Y^2,Z^2)$ is not strongly semistable on the curve given by
$Z^d=X^d+Y^d$ have density at least $(2\soger -2)/2 \soger =1-
1/\soger$.
\end{proposition}
\proof
We will show that for every remainder $r \neq 1,-1 \modu d$
there exists a power such that $r^{e} \in \rangem =\{s: d/3 < s <
d/2\}$. The residue class ring $\ZZ_d$ has $2 \soger$ units,
therefore every element has order $1$, $2$, $\soger$ or $2 \soger$.
We only have to show that $\rangem$ contains primitive remainders as
well as non-primitive remainders. Since then there exist $\varphi(2
\soger) + \varphi( \soger)= 2 \soger -2$ remainders such that some
power of them belongs to $\rangem$.

We first look for non-primitive remainders. For $d >75$ there exists
always an integer $n$ between $\sqrt{d}/\sqrt{3} <  n  <
\sqrt{d}/\sqrt{2}$, since the length of the interval is then $>1$.
Thus $n^2$ is a square in $\rangem$ and hence non-primitive. It is
also true that there exists a square in $\rangem$ for the smaller
Sophie Germain prime numbers $h=5, 11, 23, 29$.

Now to find primitive remainders note first that $d=3 \modu 4$.
Hence $-1$ is not a square in $\ZZ_d^\times$. There exist again
squares between $d/2$ and $2d/3$ (check directly for $d \leq 59$).
If $b$ is such a square, then $-b =d-b$ is a non-square inside $M$,
and so $M$ contains squares as well as non-squares (for $\soger =3
\modu 4$ one can also show by quadratic reciprocity law that $\soger
\in M$ is a non-square). \qed

\begin{remark}
If we would know that there exists infinitely many Sophie Germain
primes, then we could conclude that the density of primes for which
the bundle $\Syz(X^2,Y^2,Z^2)$ is not strongly semistable on a
Fermat curve can be arbitrarily high. The biggest Sophie Germain
prime which I have found in the literature (see
\cite{indlekoferjarai}) is $\soger =2375063906985 \cdot 2^{19380}
-1$. There should be known results in analytic number theory which
imply that the density of primes with non strongly semistable
behaviour is arbitrarily high.
\end{remark}

\begin{example}
\label{sogerfive} Let $d=11=2 \cdot 5 +1$. The set $\rangem$
consists only of $4$ and $5$. We have $2^2=4$ and $4^2= 5 \modu 11$,
hence both numbers are squares and not primitive. Hence the density
of primes $p$ such that $\Syz(X^2,Y^2,Z^2)$ is not strongly
semistable on $Z^{11}=X^{11}+Y^{11}$ for $\Char K=p$ is only $\geq
5/11$.
\end{example}

\begin{example}
Let $d=167= 2 \cdot 83 +1$. Here we have $ \rangem = \{56 \komdots
83 \}$. The numbers $s=64$ and $s=81$ are squares, hence
non-primitive elements, and $83$ is a non-square by Proposition
\ref{germainfact}. Hence the density of primes $p$ such that
$\Syz(X^2,Y^2,Z^2)$ is not strongly semistable on
$Z^{167}=X^{167}+Y^{167}$ over $\Char K=p$ is $\geq 165/167$.
\end{example}

\begin{example}
We look now at $\soger =29$, so $d= 59$. We have $\rangem =\{20
\komdots 29\}$. $2$ is a primitive element in $\ZZ_{59}$, hence
computing $2^u$, $u$ odd, (or by quadratic reciprocity law) we see
that the only primitive remainder in $\rangem$ are $23$ and $24$. So
in this range we have $8$ quadratic remainders but only $2$
non-quadratic remainders.
\end{example}

We close with an example of a degree which does not come from a
Sophie Germain prime.

\begin{example}
Let $d=31$, which does not come from a Sophie Germain prime. The
remainders $s$ for which we know that $\Syz(X^{2q},Y^{2q},Z^{2q})$
is not semistable for $q=s \modu d$, are $s \in \rangem = \{11
\komdots 15 \}$. Which remainders $p=r \modu d$ have the property
that some power $q=p^{e}=r^{e} =s =11 \komdots 15$? The number $3$
is a primitive element modulo $31$, and we have $11=3^{23}$,
$12=3^{19}$, $13=3^{11}$, $14=3^{22}$, $15=3^{21}$. So $11$, $12$
and $13$ are primitive, $14$ generates a subgroup with $15$ elements
and $15$ generates a subgroup with $10$ elements. The number of
generators of these groups are $8$, $8$ and $4$, so we have
altogether $20$ remainders for which some power fulfills the
condition in Corollary \ref{corollarynotstrong}. So the density of
primes $p$ for which $\shE_K$ is not strongly semistable in
characteristic $\Char K=p$ is at least $\geq 2/3$ (the remainders
for which we do not know the answer are $1,2,4,5,8,16,
25,27,29,30$).
\end{example}

\bibliographystyle{plain}

\bibliography{bibliothek}

\end{document}